# Observations about Leonardo's drawings for Luca Pacioli

DIRK HUYLEBROUCK

Department for Architecture, KULeuven, Brussels, Belgium.

*Three versions of Luca Pacioli's 'De Divina Proportione' remain: a manuscript held in Milan, another in Geneva and a printed version edited in Venice. A recent book, 'Antologia della Divina Proporzione', has all three in one volume, allowing an easy comparison of the different versions. The present paper proposes some observations about these drawings, generally said to be of Leonardo da Vinci's hand.*

## Three copies

Three versions of Luca Pacioli's *'De divina proportione'* remain: there are the two manuscripts, one from the Bibliothèque de Genève and one from the Biblioteca Ambrosiana in Milan, written between 1496 and 1498, and there is the wider spread printed version from Venice, which appeared only in 1509, that is, about ten years later. A recent book, 'Antologia della Divina Proporzione' unites images from all three editions and allows comparing them easily by presenting similar illustrations on opposing pages (see [1]). Not surprisingly, the colour and printed images are each other's reflections about a vertical axis (in most cases) and thus we will represent them here often with an additional reflection about the vertical axis, to make the comparison easier.

The present paper focuses on observations based on the comparison of the different version. This is not straightforward: for instance, the colour plates are not always as accurate as the black and white ones, or vice versa. Still, the present approach is rather modest, as we avoid interpretations and historic questions, which we leave to art historians. To mathematicians, the proposed comparisons are interesting and intriguing enough and one can only wonder why they were not made before. For instance, in the works by Kim Veltman (see [9]) and Ladislao Reti (see [8]) many examples are given, from the different versions, and although there sometimes is some comparison of the geometric representations, no explicit or detailed observations appear in these yet voluminous works. Yet, Reti for instance does not hesitate to systematically point out some problems in Leonardo's work when it comes to squaring the circle or drawing ballistic trajectories.

Of course, there are other regular polyhedra in some of Leonardo's codices, and here too careful observation is sometimes interesting. In the Codex Arundel for instance figures an icosahedron that seems to have 13 vertices (see [6]). However, these codices were draft writings, and no official publications as Pacioli's Divina was one. Sometimes, they were rough copies, as on folio 707 of the Codex Atlanticus, where a drawing of an icosidodecahedron seems to be glued on the page. It is identical to the drawing Leonardo made for Pacioli, but turned under an angle. It could be interesting to try to find out more about these other drawings (was the glued drawing one of the originals Leonardo made for Pacioli?), but they are not the topic of the present paper.

## Front views

Of most polyhedra Leonardo drew for Pacioli, he made four versions: a 'planar' version, showing the faces, a 'vacuus' version, showing the edges, and two 'stellar' versions, again in planar and vacuus version, where pyramids where drawn on each face. Thus, one would expect similarities as given for the cube, where the planar and the skeleton correspond perfectly, in all three editions, the two manuscripts and the printed version (see fig. 1).

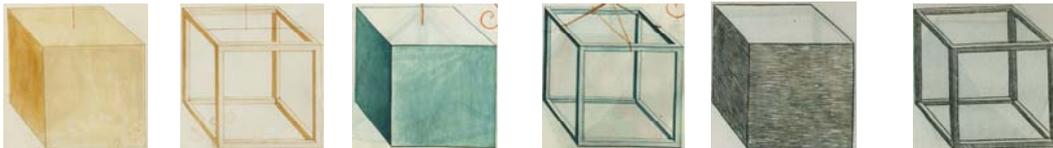

*Fig. 1: The Geneva, Milan and (mirrored) printed versions of the cube, in planar version and skeleton version.*

However, this similarity is not always followed. An art historian who already pointed this out in one case is Noémie Etienne (see [2]): 'The solid prism with triangular basis (the number XLIII in the Geneva manuscript) is represented in a different way in the three works: in the Geneva copy, the angle is in front of the spectator, while its representation with edges shows a flat face to the spectator. In the Milan manuscript, the same illustration is represented twice from the same viewpoint: it is seen from the back, in the solid representation as well as in the edge form. Finally, in the 1509 edition, the illustration is presented in a similar way as in the Geneva manuscript, though the solid was turned upside down' (see fig. 2).

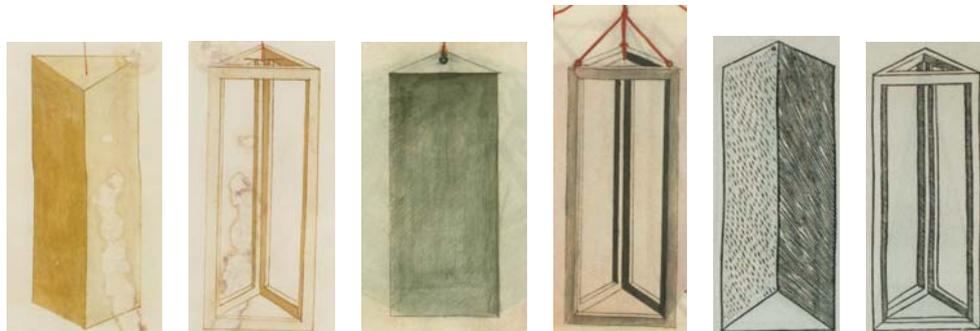

*Fig. 2: Geneva, Milan and (mirrored) printed versions of the triangular prism.*

The latter conclusion, that the planar triangular prism was turned upside down is speculation: the tetrahedron has a similar representation in the printed version and clearly it was not reversed. Also, it remains surprising this observation about the planar tetrahedron was not made, while the previous one was (see fig. 3). Note even the Geneva and the Milan planar and skeleton versions do not really correspond for the tetrahedron.

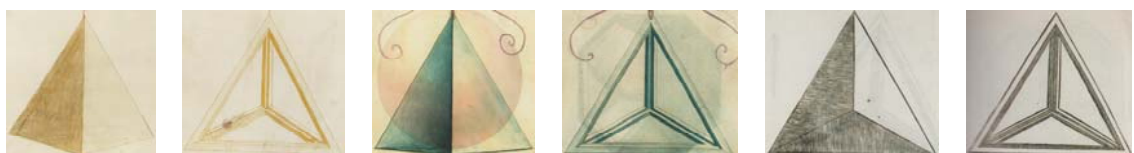

*Fig. 3: The Geneva, Milan and (mirrored) printed versions of the tetrahedron.*

For the 'lateral triangular pyramid' the same front-behind switches are made (see fig. 5).

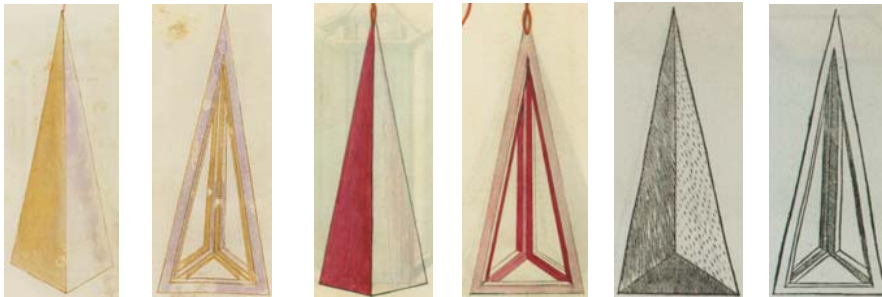

*Fig. 4: The Geneva, Milan and (mirrored) printed versions of a lateral triangular pyramid.*

In the case of the 'lateral triangular non-equilateral pyramid' the situation is even more noteworthy, since there is an additional left-right switch for the back vertex, in the printed version (see fig. 5).

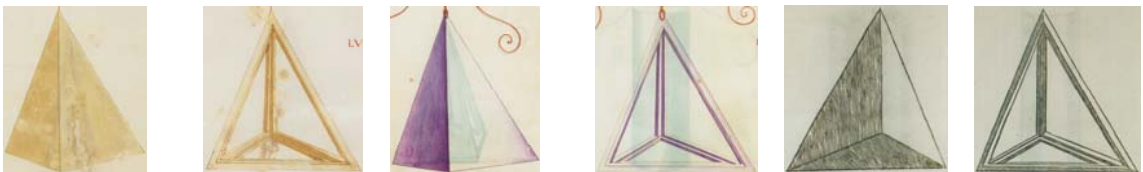

*Fig. 5: The Geneva, Milan and (mirrored) printed versions of a lateral triangular non-equilateral pyramid.*

On the other hand, the versions of the 'lateral square pyramid' do correspond in the Geneva and Milan copy, but there now there is a symmetry flip in the printed version (see fig. 6).

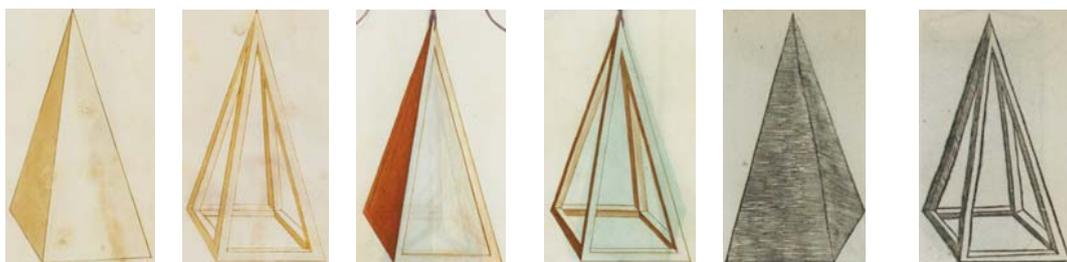

*Fig. 6: The Geneva, Milan and (mirrored) printed versions of a lateral square pyramid.*

The explanations given by Pacioli in his text compare the prismatic solids to the top of the roof of a house constituted by 4 faces (see section LIX of Pacioli's printed text). The drawing of a triangular prism given in the margin of that paragraph indeed looks like a top view on a roof (see fig. 7).

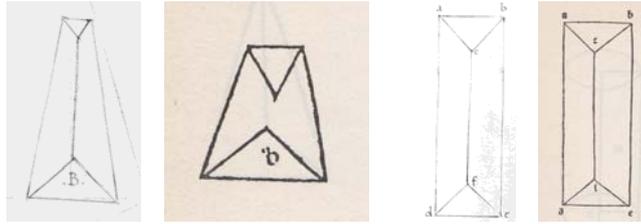

*Fig. 7: The drawings in the margins of the text in De Divina Proportione: the truncated triangular pyramid in the Geneva (left) and printed version (middle left) and triangular prism in the Geneva (middle right) and printed version (right).*

## Shadows

Shadows on polyhedra are very accurately represented in the Milan version, while in the Geneva version they are vaguer. The reader is invited to observe the printed version, for instance for the exacedron abscisus elevatus (fig. 8). These are not the only such cases, but they are rather striking.

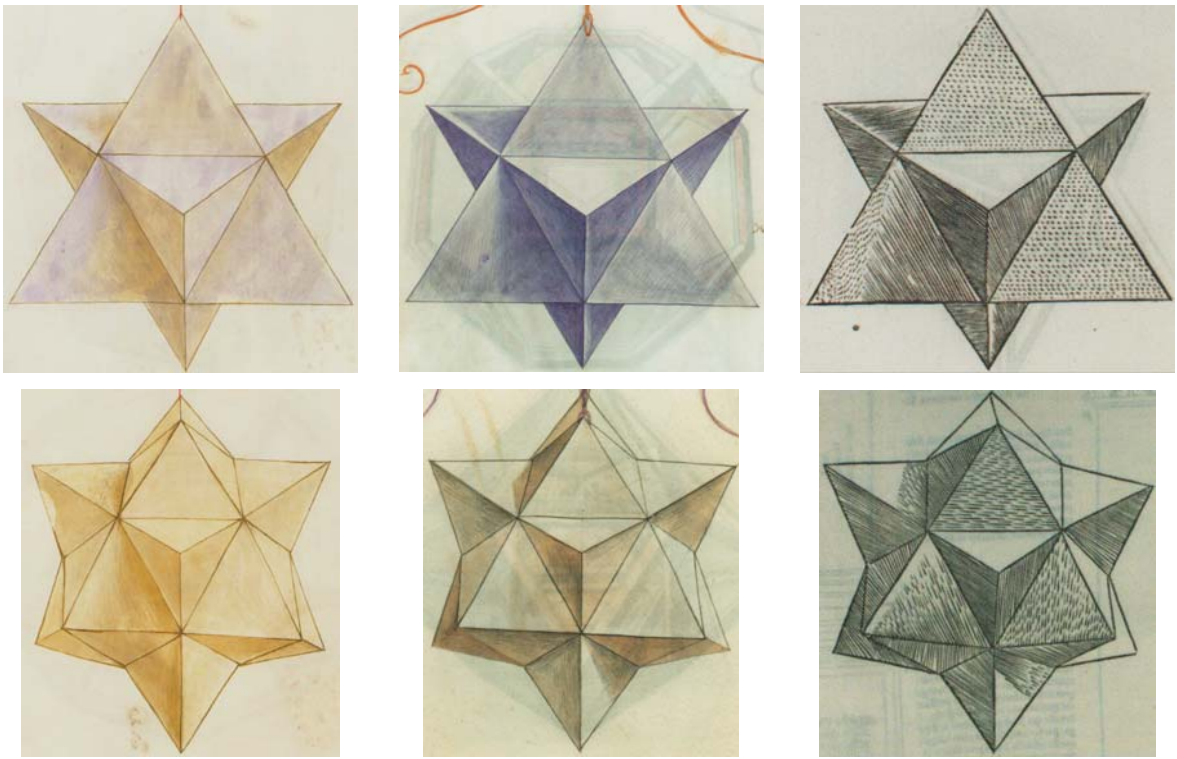

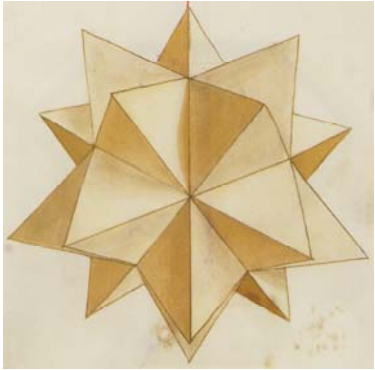 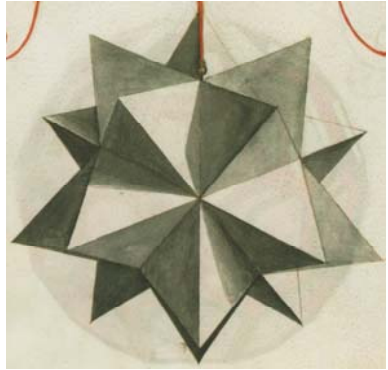 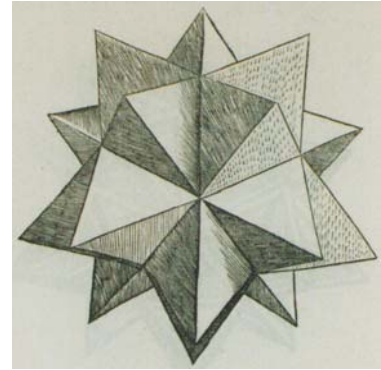

*Vague shadows*          *Perfect shadows*          *Wrong shadows*

*Fig. 8: Compare the shadows on the Geneva, Milan and (mirrored) printed versions of the 'octohedron elevatus' (above), the 'exacedron abscisus elevatus' (middle) and 'ycocedron elevatus' (below).*

## Faces and ropes

One may have the impression the Geneva version is the most elegant one of the three versions. After all, it was the copy Pacioli offered to the Duke of Milan. However, some observations can be made about this version too: some faces of polyhedra in skeleton versions are planar.

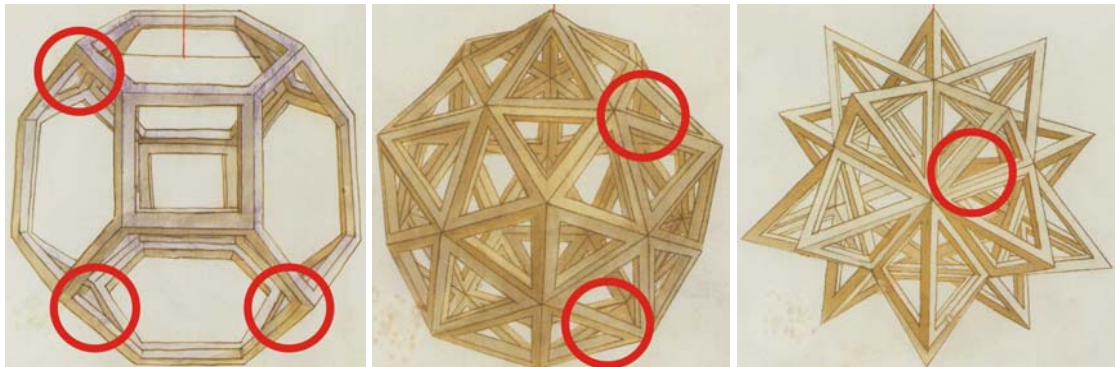

*Fig. 9: Some faces of polyhedra in the Geneva skeleton versions are 'planar'.*

Moreover, the ropes at which the polyhedra hung provide another topic of discussion. Some suspensions in the Geneva version seem merely decorative, while they look rather realistic in the Milan version. The printed version shows no ropes at all.

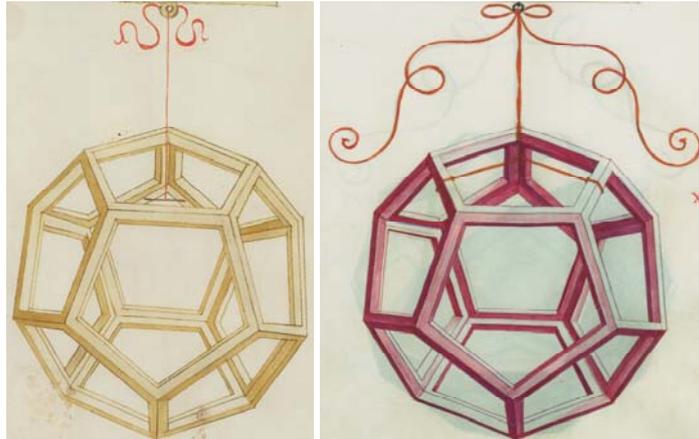
*Fig. 10: The ropes in the Geneva version (left), and in the Milan version (right).*

## Printed version: rotations and tangents

Leonardo did not appreciate low quality woodcuts (see [5]). On Folio 139v of his 'Anatomical Studies' he added the note, dated at about 1510 (that is, after the printing of 'De Divina Proportione'): "I beg the one coming after me not to be guided by thrift and use a wood cut". This low quality is especially illustrated for the sphere, which hardly looks like one in the woodcut version. Also, when drawing a cone, the base circle should be an ellipse, and the lines to the apex should be tangent lines to this ellipse. This is hard to see in the Geneva and Milan version and doubtful in the printed version: the straight lines down from the apex to the base circle are almost orthogonal to it. Note that in the drawings in the margins, for instance for the cylinder, the circles of the bottom and top surfaces are no ellipses, in all three versions.

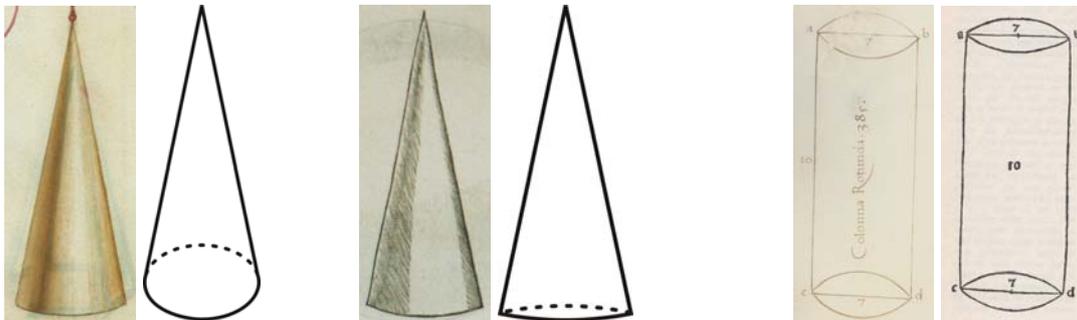
*Fig. 11: The cone in the Milan version (left) and in the (mirrored) printed version (middle; the drawings are exaggerated); cylinders (right).*

Also noteworthy is that the printed version of the icosidodecahedron is rotated over 72°, and mirrored, when compared to the Geneva or the Milan version.

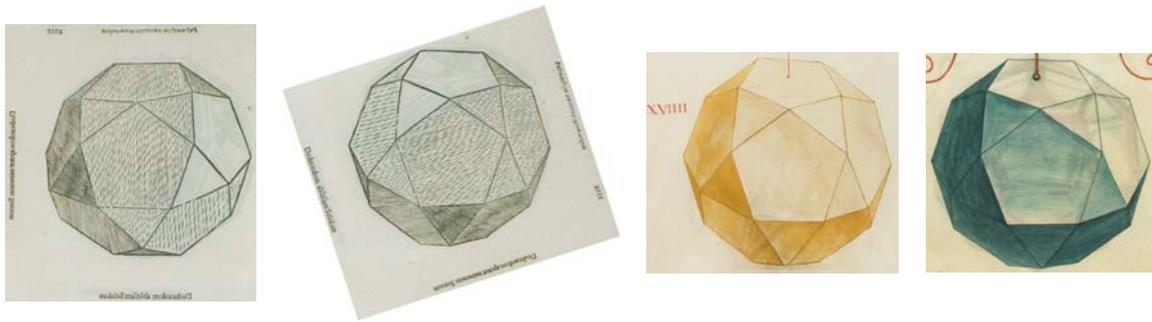

*Fig. 12: The (mirrored) printed version of the icosidodecahedron (left) and the (not mirrored) rotated version (middle left), compared to the Geneva version (middle right) and the Milan version (right).*

A similar observation holds for the skeleton version of the stellar dodecahedron, as noticed by Jos Janssen. However, just as in the previous case, there is an additional reflection of the image.

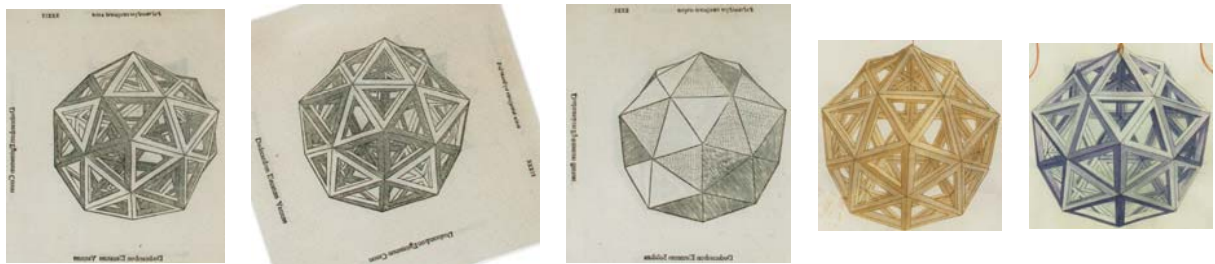

*Fig. 13: The (mirrored) printed version of the skeleton stellar dodecahedron (left) and the (not mirrored) rotated version (middle left), compared to the (mirrored) planar version (middle) and the Geneva and Milan version (right).*

## Geneva: narrower drawings

Below is the list of the polyhedra from the Divina mentioned in the '*Antologia della Divina Proporzione*' (see table 1). There are a few more, since there is a difference between the Geneva and the Milan version. There are some minor differences in the names too: in the Geneva version, the label says the 'exacedron' is also called 'cube', while the 'Columna laterata triangula solida' and 'vacua' both get the additional description 'seu corpus seratile'. Yet, more remarkable is that the proportion of the width divided by the length of the drawings is always smaller in the Geneva version. The drawings in the Geneva version look narrower, and this can indeed be observed by the naked eye in some cases (the measurements and numbers were only made for an additional verification). The printed versions seem often smaller too, with respect to the Milan version, but here this verification was only made for the sake of completeness, as the printed versions are rather small and thus the measurement do not make that much sense in comparison to the other measurements.

| Name volume | Width /Length Geneva | Width /Length Milan | Width /Length Printed | Difference Milan - Geneva | Difference Milan - printed | Difference Geneva - printed |
|---|---|---|---|---|---|---|
| Tetracedron planus solidus | 1.02 | 1.09 | 1.17 | 0.07 | -0.08 | -0.15 |
| Tetracedron planus vacuus | 1.13 | 1.19 | 1.15 | 0.06 | 0.04 | -0.02 |
| Exahedron planus solidus | 0.94 | 1 | 0.98 | 0.06 | 0.02 | -0.04 |
| Exahedron planus vacuus | 0.97 | 1 | 0.98 | 0.03 | 0.02 | -0.01 |

| Name | | | | | | |
|---|---|---|---|---|---|---|
| Tetracedron abscisus solidus | 0.95 | 1.02 | 0.96 | 0.07 | 0.06 | -0.01 |
| Tetracedron abscisus vacuus | 0.98 | 1.03 | 0.96 | 0.05 | 0.07 | 0.02 |
| Ycocedron elevatus solidus | 0.99 | 1.03 | 1 | 0.04 | 0.03 | -0.01 |
| Ycocedron elevatus vacuus | 1.02 | 1.02 | 1.02 | 0 | 0 | 0 |
| Exacedron abscisus solidus | 1.14 | 1.17 | 0.96 | 0.03 | 0.21 | 0.18 |
| Exacedron abscisus vacuus | 1.13 | 1.18 | 0.95 | 0.05 | 0.23 | 0.18 |
| Duodecedron abscisus solidus | 1.01 | 1.04 | 1.04 | 0.03 | 0 | -0.03 |
| Duodecedron abscisus vacuus | 1 | 1.05 | 1.04 | 0.05 | 0.01 | -0.04 |
| Octocedron abscisus solidus | 0.99 | 1.02 | 1 | 0.03 | 0.02 | -0.01 |
| Octocedron abscisus vacuus | 0.98 | 1.01 | 1.01 | 0.03 | 0 | -0.03 |
| Octocedron planus solidus | 0.83 | 0.97 | 0.96 | 0.14 | 0.01 | -0.13 |
| Octocedron planus vacuus | 0.93 | 0.97 | 0.93 | 0.04 | 0.04 | 0 |
| Vigintisex basium planus solidus | 0.8 | 0.92 | 0.88 | 0.12 | 0.04 | -0.08 |
| Vigintisex basium planus vacuus | 0.87 | 0.92 | 0.89 | 0.05 | 0.03 | -0.02 |
| Ycocedron planus solidus | 0.9 | 0.93 | 0.91 | 0.03 | 0.02 | -0.01 |
| Ycocedron planus vacuus | 0.92 | 0.93 | 0.95 | 0.01 | -0.02 | -0.03 |
| Exacedron elevatus solidus | 1.1 | 1.11 | 1.13 | 0.01 | -0.02 | -0.03 |
| Exacedron elevatus vacuus | 1.1 | 1.12 | 1.12 | 0.02 | 0 | -0.02 |
| Duodecedron elevatus solidus | 0.96 | 0.96 | 1.04 | 0 | -0.08 | -0.08 |
| Duodecedron elevatus vacuus | 0.96 | 0.98 | 1.04 | 0.02 | -0.06 | -0.08 |
| Octocedron elevatus solidus | 0.96 | 0.99 | 1.02 | 0.03 | -0.03 | -0.06 |
| Octocedron elevatus vacuus | 0.99 | 0.99 | 0.88 | 0 | 0.11 | 0.11 |
| Duodecedron planus solidus | 0.96 | 0.96 | 0.94 | 0 | 0.02 | 0.02 |
| Duodecedron planus vacuus | 0.93 | 0.95 | 0.96 | 0.02 | -0.01 | -0.03 |
| Exacedron abscisus elevatus solidus | 0.83 | 0.85 | 0.91 | 0.02 | -0.06 | -0.08 |
| Exacedron abscisus elevatus vacuus | 0.82 | 0.84 | 0.83 | 0.02 | 0.01 | -0.01 |
| Vigintisex basium elevatus solidus | 0.98 | 1.03 | 0.99 | 0.05 | 0.04 | -0.01 |
| Vigintisex basium elevatus vacuus | 0.99 | 1.01 | 1 | 0.02 | 0.01 | -0.01 |
| Tetracedron elevatus solidus | 1.13 | 1.19 | 1.17 | 0.06 | 0.02 | -0.04 |
| Tetracedron elevatus vacuus | 1.17 | 1.19 | 1.2 | 0.02 | -0.01 | -0.03 |
| Ycocedron abscisus solidus | 1 | 1.01 | 1.02 | 0.01 | -0.01 | -0.02 |
| Ycocedron abscisus vacuus | 1.02 | 1.02 | 1 | 0 | 0.02 | 0.02 |
| Duodecedron abscisus elevatus solidus | 0.96 | absent | 0.92 | | | 0.04 |
| Duodecedron abscisus elevatus vacuus | 0.93 | 0.95 | 0.93 | 0.02 | 0.02 | 0 |
| Septuaginta duarum basium solidum | 0.99 | 1.02 | 1 | 0.03 | 0.02 | -0.01 |
| Septuaginta duarum basium vacuum | 1 | 1.02 | 1 | 0.02 | 0.02 | 0 |
| Sphera solida | 1 | 1 | 1 | 0 | 0 | 0 |
| Columna laterata triangula solida seu corpus seratile and Columna laterata triangula solida | 0.39 | 0.39 | 0.38 | 0 | 0.01 | 0.01 |
| Columna laterata triangula vacua seu corpus seratile and Columna laterata triangula vacua | 0.38 | 0.39 | 0.38 | 0.01 | 0.01 | 0 |
| Columna laterata quadrangula solida | 0.5 | 0.51 | 0.89 | 0.01 | -0.38 | -0.39 |
| Columna laterata quadrangula vacua | 0.5 | 0.5 | 0.52 | 0 | -0.02 | -0.02 |
| Columna laterata pentagona solida | 0.46 | 0.46 | 0.47 | 0 | -0.01 | -0.01 |
| Columna laterata pentagona vacua | 0.46 | 0.47 | 0.46 | 0.01 | 0.01 | 0 |
| Columna laterata exagona solida | 0.48 | 0.48 | 0.48 | 0 | 0 | 0 |
| Columna laterata exagona vacua | 0.49 | 0.49 | 0.49 | 0 | 0 | 0 |
| Piramis laterata triangula solida | 0.35 | 0.39 | 0.4 | 0.04 | -0.01 | -0.05 |
| Piramis laterata triangula vacua | 0.39 | 0.39 | 0.4 | 0 | -0.01 | -0.01 |
| Pi(y)ramis laterata triangula inequilatera solida | 0.86 | 0.9 | 0.98 | 0.04 | -0.08 | -0.12 |
| Pyramis laterata triangula inequilatera vacua | 0.96 | 0.98 | 0.98 | 0.02 | 0 | -0.02 |
| Pyramis laterata pentagona solida | 0.49 | 0.5 | 0.49 | 0.01 | 0.01 | 0 |
| Pyramis laterata pentagona vacua | 0.48 | 0.49 | 0.49 | 0.01 | 0 | -0.01 |
| Columna rotunda solida | 0.32 | 0.34 | 0.34 | 0.02 | 0 | -0.02 |
| Pyramis rotunda solida | 0.36 | 0.37 | 0.39 | 0.01 | -0.02 | -0.03 |
| Piramis laterata quadrangula solida | 0.54 | 0.55 | 0.58 | 0.01 | -0.03 | -0.04 |
| Piramis laterata quadrangula vacua | 0.54 | 0.55 | 0.58 | 0.01 | -0.03 | -0.04 |

*Table 1: List of the polyhedra represented in [1] and their width over length quotient.*

The drawings are systematically narrower in the Geneva version than in the Milan version (see fig. 14).

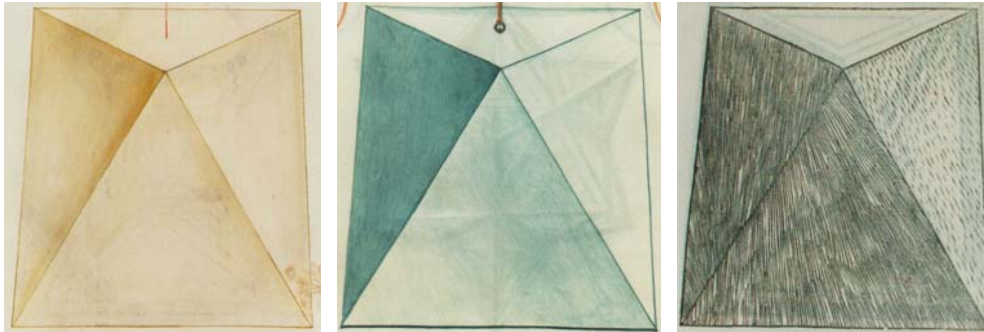

*Fig. 14: Reduced to an identical height, the drawings in the Milan manuscript (middle) are larger than those in the Geneva manuscript (right) and also often larger than those in the printed version (right, mirrored image).*

## Observing the edges

It is interesting too to observe the edges in different versions of polyhedra. For instance, a more detailed look on the stellar icosidodecahedron or 'duodecedron abscisus elevatus vacuus' shows problems on the lower pyramid below the central pentagonal pyramid in the three versions. The Milan version seems to miss two more edges, while the printed version lacks one edge (though it could be a matter of interpretation).

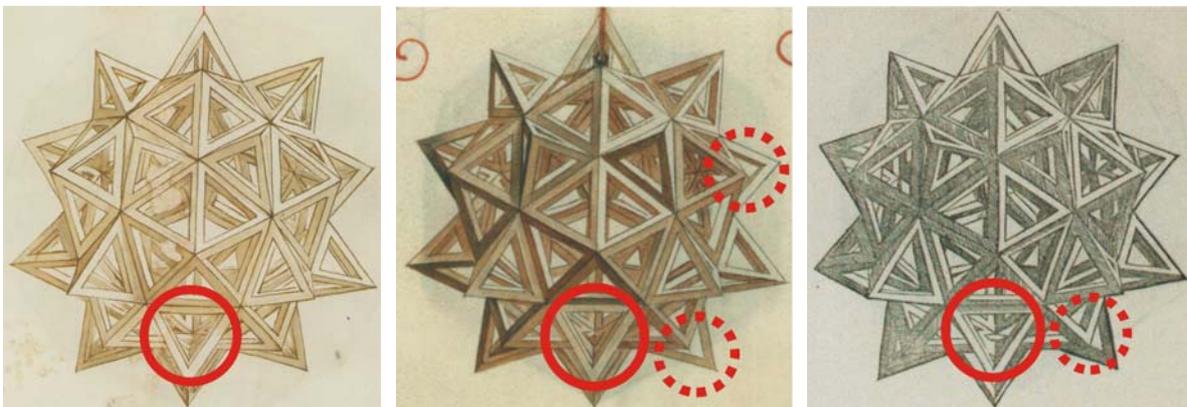

*Fig. 15: The Geneva, Milan and (mirrored) printed version of the stellar icosidodecahedron and some indications.*

A similar problem was noticed by chemist Jos Janssen (The Netherlands) while discussing some topics of the present paper. The stellar version of the cuboctahedron looks like an Escher or Penrose drawing: the edges coming from the triangle in the back should be behind all other edges, but suddenly lay in front when they reach the outer top vertices of the triangular pyramids they are composing. Also, observe the top (horizontal) edge of the square on top of the cuboctahedron is in the stellar version, in all three versions of the De Divina Proportione.

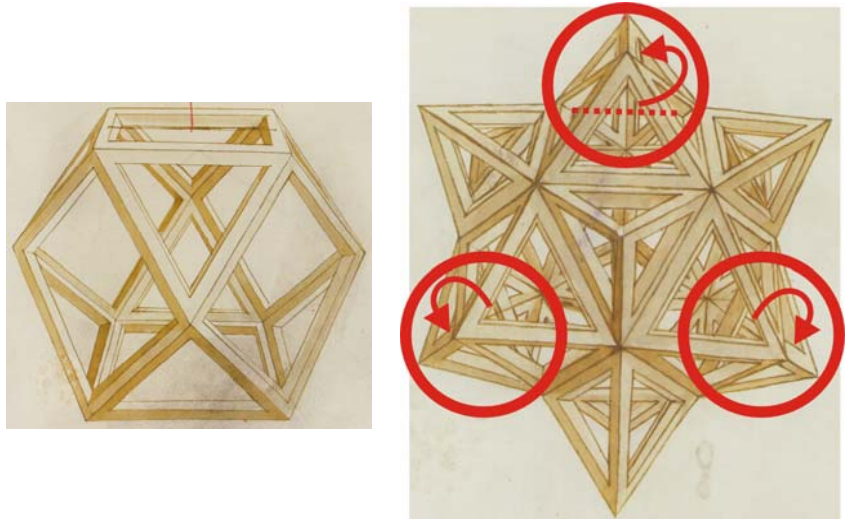

*Fig. 16: The cuboctahedron and the stellar cuboctahedron in the Geneva version.*

The observations about the stellar icosidodecahedron and the stellar cuboctahedron remind a problem noticed by Rinus Roelofs about the stellar rhombohedric cuboctahedron (see [3], [4]). It was called a 'geometric error' and many Leonardo enthusiasts, sometimes motivated by his alleged infallible status, went at length in their efforts to explain it: it would have been a stellar pseudo-rhombohedric cuboctahedron; it would have been a riddle; or even a mystical message.

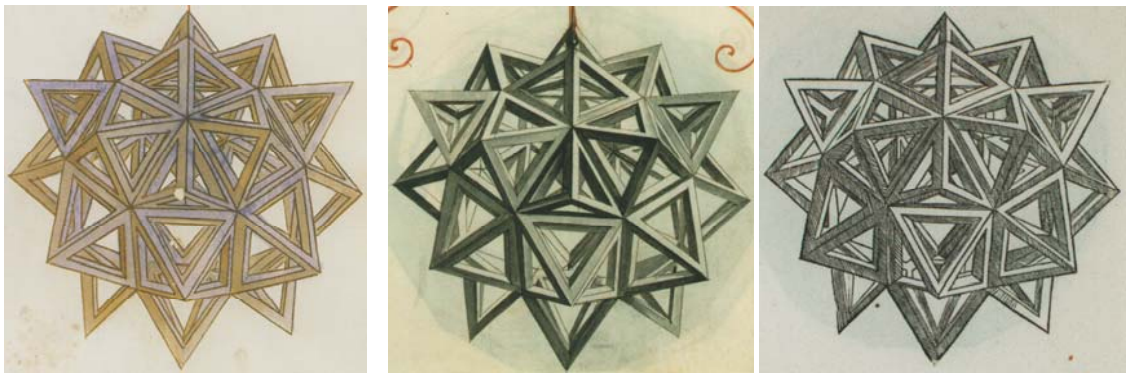

*Fig. 17: The stellar rhombohedric cuboctahedron in the three versions; observe the lower pyramid (right: mirrored image).*

However, we learned from the observations made in the present paper, that the inaccuracy in this drawing is but one of many, and that it indeed appears to be nothing more than that, that is, an inaccuracy. It was the first time polyhedra were so clearly drawn in colour, and with a high suggestive 3D-effect. Today, this spatial appearance is easily obtained on computer, and it was this way expert computer draftsman Rinus Roelofs noticed the problem with the rhombohedric cuboctahedron. However, after all, and viewed in his time, it should perhaps not be seen as a 'geometric error' illustrating some 'mathematical mistake', but as a simple 'drawing inaccuracy' – which, perhaps was not even Leonardo's but a copyists' or an engraver's, but that is yet another discussion.